\documentclass[11pt]{article}
\everymath{\displaystyle}
\usepackage{amssymb}
\usepackage{amsmath}
\newtheorem{prethm}{{\bf Theorem}}

\newenvironment{thm}{\begin{prethm}{\hspace{-0.5
               em}{\bf .}}}{\end{prethm}}
\newtheorem{prelemma}{{\bf Lemma}}

\newenvironment{lemma}{\begin{prelemma}{\hspace{-0.5
               em}{\bf .}}}{\end{prelemma}}
\newtheorem{preex}{{\bf Example}}

\newenvironment{ex}{\begin{preex}{\hspace{-0.5
               em}{\bf .}}}{\end{preex}}
\newtheorem{preprop}{{\bf Proposition}}

\newenvironment{prop}{\begin{preprop}{\hspace{-0.5em}{\bf .}}}{\end{preprop}}
\newtheorem{precor}{{\bf Corollary}}

\newenvironment{cor}{\begin{precor}{\hspace{-0.5
               em}{\bf .}}}{\end{precor}}
\newtheorem{preremark}{{\bf Remark}}

\newenvironment{remark}{\begin{preremark}{\hspace{-0.5
               em}{\bf.}}}{\end{preremark}}
\newtheorem{preprob}{{\bf Problem}}

\newtheorem{predefin}{{\bf Definition}}

\newtheorem{preconj}{{\bf Conjecture}}

\newtheorem{preprobb}{{\bf Problem}}

\newtheorem{prelem}{{\bf Theorem}}

\newenvironment{proof}{{\bf Proof.}\rm }{\hfill{$\Box$}}

\newtheorem{presolution}{{\bf Solution.}}

\def\newpic#1{}
\def\qed{\ifhmode\unskip\nobreak\fi\quad\ifmmode\Box\else$\Box$\fi}
\title{\vspace{-1cm}\Large\bf\noindent Lower bounds for independence and $k$-independence number of graphs using the concept of degenerate degrees}

\author{\large\bf Manouchehr Zaker\footnote{mzaker@iasbs.ac.ir}
\vspace{5mm}\\
    Department of Mathematics,\\
     Institute for Advanced Studies in Basic Sciences,\\
    Zanjan 45137-66731, Iran\\
  }
    \date{}
\begin{document}
\maketitle
\begin{abstract}
\noindent Let $G$ be a graph and $v$ any vertex of $G$. We define the degenerate degree of $v$, denoted by $\zeta(v)$ as $\zeta(v)={\max}_{H: v\in H}~\delta(H)$, where the maximum is taken over all subgraphs of $G$ containing the vertex $v$. We show that the degenerate degree sequence of any  graph can be determined by an efficient algorithm. A $k$-independent set in $G$ is any set $S$ of vertices such that $\Delta(G[S])\leq k$. The largest cardinality of any $k$-independent set is denoted by $\alpha_k(G)$. For $k\in \{1, 2, 3\}$, we prove that $\alpha_{k-1}(G)\geq {\sum}_{v\in G} \min \{1, 1/(\zeta(v)+(1/k))\}$. Using the concept of cheap vertices we strengthen our bound for the independence number. The resulting lower bounds improve greatly the famous Caro-Wei bound and also the best known bounds for $\alpha_1(G)$ and $\alpha_2(G)$ for some families of graphs. We show that the equality in our bound for independence number happens for a large class of graphs. Our bounds are achieved by Cheap-Greedy algorithms for $\alpha_k(G)$ which are designed by the concept of cheap sets. At the end, a bound for $\alpha_k(G)$ is presented, where $G$ is a forest and $k$ an arbitrary non-negative integer.
\end{abstract}


\section{Introduction and related works}

\noindent Degeneracy of graphs is a very useful and widely studied concept in graph theory which was firstly introduced by Lick and White in \cite{LW}. For any two graphs $H$ and $G$, by $H\subseteq G$ we mean $H$ is a subgraph of $G$. If $S$ is any subset of vertices in $G$, we denote by $G[S]$ the subgraph of $G$ induced by the elements of $S$. Let $G$ be any graph and $v$ any vertex of $G$. Denote the degree of $v$ in $G$ by $deg_G(v)$. A graph $G$ is $k$-degenerate if its vertices can be ordered as $v_1, v_2, \ldots, v_n$ such that for any $i$, $deg_{G_i}v_i \leq k$, where $G_i=G[v_1, \ldots, v_i]$. It was also known that the minimum $k$ such that $G$ is $k$-degenerate equals to ${\max}_{H:H\subseteq G}~\delta(H)$, where $\delta(H)$ denotes the minimum degree of $H$. The latter value is called the degeneracy of $G$. The terminology {\it degenerate graph} sounds somewhat strange. In a personal communication with A. White, an author of \cite{LW} mentioned the following comment concerning the terminology degenerate graph. {\it ``My thinking was that: ``degenerate" means deteriorated, that is passed from a higher quality to a lesser one. Thus a $k$-degenerate graph deteriorates to the empty graph by a sequence of deleting vertices of degree at most $k$ (in each step in the sequence)"} \cite{Wh}. Some problems which are NP-hard for general graphs have polynomial time solutions for degenerate graphs (see e.g. \cite{AG}). Degeneracy is useful idea in chromatic and extremal graph theory and caused many interesting problems and topics (see e.g. \cite{JT}). The complement (with respect to the vertex set) of degenerate subgraphs possess another interesting graph theoretical property i.e. the so-called dynamic monopolies \cite{Z}.

\noindent In this paper we introduce the {\it degenerate degree} of vertices as follows. Let $G$ be a graph and $v$ any vertex of $G$. Define the degenerate degree of $v$ in $G$ as ${\max}_{H: v\in H}~\delta(H)$, where $\delta(H)$ stands for the minimum degree of $H$ and the maximum is taken over all subgraphs of $G$ which include the vertex $v$. We denote this value by $\zeta_G(v)$ (or simply by $\zeta(v)$ when the underlying graph is known and fixed). Note that if $v$ is an isolated vertex in a graph, then $\zeta(v)=0$. It is clear by the definition that for any vertex $\delta(G)\leq \zeta(v)\leq deg_G(v)$. Denote the neighborhood set of any vertex $v$ in $G$ by $N(v)$. For any subset $S$ of the vertices of $G$, by $N(S)$ and $N[S]$ we mean respectively, ${\cup}_{v\in S} N(v)$ and $N(S)\cup S$. Proposition \ref{algor} provides an efficient algorithm which outputs the degenerate degree sequence of graphs. In the following proposition we use the concept of {\it smallest-last order}. An ordering $v_1, v_2, \ldots, v_n$ of the vertices of a graph $G$ is said to be a smallest-last order if $v_i$ is a vertex with the minimum degree in the graph $G\setminus \{v_1, \ldots, v_{i-1}\}$, for any $i\in \{1, \ldots, n\}$.
It is a well-known fact that a smallest-last vertex ordering for a graph on the vertex set $V$ and the edge set $E$ can be obtained in linear time complexity $\mathcal{O}(|V|+|E|)$ using space complexity $2|E|+\mathcal{O}(|V|)$ (see e.g. \cite{MB}).

\begin{prop}
There exists an $\mathcal{O}(|E|+|V|)$ algorithm which outputs the sequence $(\zeta(v))_{v\in G}$ for any input graph $G$ with the vertex set $V$ and the edge set $E$.\label{algor}
\end{prop}

\noindent \begin{proof}
Let $v_1, v_2, \ldots, v_n$ be a smallest-last vertex ordering in $G$. Set for simplicity $\tilde{d}_i=deg_{G[v_i, \ldots, v_n]}(v_i)$. We show that
$$\zeta(v_i) =\max_{j=1}^i \tilde{d}_j.$$
\noindent Let $j$ be any value with $1\leq j\leq i$. Then $v_i\in G[v_j, \ldots, v_n]$ and by the definition $\zeta(v_i) \geq \delta(G[v_j, \ldots, v_n])=\tilde{d}_j$. This shows that $\zeta(v_i) \geq {\max}_{j\in \{1, \ldots, i\}}~\tilde{d}_j$. Now let $H$ be a subgraph of $G$ such that $v_i\in H$ and $\zeta(v_i)=\delta(H)$. Let $k\in \{1, \ldots, i\}$ be such that $v_k\in H$ but $V(H)\cap \{v_1, \ldots, v_{k-1}\}=\varnothing$. Then $H\subseteq G[v_k, \ldots, v_n]$ and $deg_H(v_k)\leq \tilde{d}_k$. It follows that $\delta(H)\leq {\max}_{j\in \{1, \ldots, i\}}~\tilde{d}_j$. This proves $\zeta(v_i) = {\max}_{j\in \{1, \ldots, i\}}~\tilde{d}_j$, for each $i=1, \ldots, n$.

\noindent Now we provide an algorithm which outputs the sequence $(\zeta(v))_{v\in G}$ for any graph $G$. As we mentioned earlier, a smallest-last vertex ordering $v_1, v_2, \ldots, v_n$ can be obtained in $\mathcal{O}(|V|+|E|)$ time steps. Therefore, the values $\tilde{d}_i$ can be determined with the same time complexity. The equality $\zeta(v_i) = {\max}_{j\in \{1, \ldots, i\}}~\tilde{d}_j$ shows that
for each $i=0, 1, \ldots, n-1$, $\zeta(v_{i+1}) = \max \{\zeta(v_i),~\tilde{d}_{i+1}\}$. This means that once $\zeta(v_{i})$ is determined, $\zeta(v_{i+1})$ can also be determined from $\zeta(v_i)$ by only one comparison. We conclude that after determining the sequence $(\tilde{d}_i)_{i=1}^n$, the sequence $(\zeta(v))_{v\in G}$ can be computed in $\mathcal{O}(|V|)$ time steps. The total time complexity is $\mathcal{O}(|E|+|V|)$. This completes the proof.
\end{proof}

\noindent We call a graph $G$, $\zeta$-regular if $\zeta(u)=\zeta(v)$ for any two vertices $u, v$ of $G$. The following lemma concerns $\zeta$-regular graphs.

\begin{lemma}
A graph $G$ is $\zeta$-regular if and only if $\max_{K\subseteq G} \delta(K) = \delta(G)$.
\end{lemma}

\noindent \begin{proof}
Assume first that for any two vertices $u$ and $v$, $\zeta(u)=\zeta(v)$. Let ${\max}_{K\subseteq G} \delta(K) = \delta(H)$. It is enough to show that $\delta(H)\leq \delta(G)$. Let $u_0$ and $v_0$ be such that $deg_G(u_0)=\delta(G)$ and $deg_H(v_0)=\delta(H)$. We have $\delta(H)\leq \zeta(v_0) = \zeta(u_0) \leq deg_G(u_0) = \delta(G)$. Conversely, assume that $\delta(H) = \delta(G)$. Let $v$ be any vertex of $G$. We have $\delta(G)\leq \zeta(v)\leq {\max}_{K\subseteq G} \delta(K) = \delta(H) =\delta(G)$. It follows that for any vertex $\zeta(v)=\delta(G)$, as desired.
\end{proof}

\noindent Throughout the paper we use some special vertices which we call {\it cheap vertices}. A vertex $u$ is said to be a cheap vertex if $\zeta(u) = deg_G(u)$ and $\zeta(u)= {\min}_{v\in N[u]}\zeta(v)$. Note that we can replace the latter condition by $\zeta(u) \leq {\min}_{v\in N(u)}\zeta(v)$. It is easy to observe that any vertex with minimum degree in $G$ is a cheap vertex of $G$. The cheap vertices are not limited to vertices with minimum degree. For example consider a graph $G$ consisting of a triangle with a pendant vertex (i.e. a vertex of degree one in $G$). Then the two vertices of degree two in $G$ are cheap and the minimum degree of the graph is one. We see in this paper the importance of cheap vertices including those which have not the minimum degree. The following proposition shows that the set of cheap vertices in any given graph can be found by an efficient algorithm.

\begin{prop}
Let $G$ be any given graph on the vertex and edge set $V$ and $E$, respectively. The set of all cheap vertices in $G$ can be obtained by an efficient algorithm with time complexity $\mathcal{O}(|E|+|V|)$.\label{output-cheap}
\end{prop}

\noindent \begin{proof}
By Proposition \ref{algor}, there exists an $\mathcal{O}(|E|+|V|)$ algorithm which determines the degenerate degree sequence of $G$. Using this algorithm we obtain the set of vertices $v$ satisfying $deg_G(v)=\zeta(v)$. For each such vertex $v$ we check wether $\zeta(v)= {\min}_{w\in N[v]}\zeta(w)$. This can be done in $deg_G(v)$ time steps for each vertex $v$ and at most ${\sum}_{v\in G} deg_G(v)=2|E|$ time steps for all other vertices. Therefore the total time complexity is $\mathcal{O}(|E|+|V|)$.
\end{proof}

\noindent Let $k$ be any non-negative integer. In this paper we deal with the $k$-independence number of graphs. A subset $S$ of the vertices of $G$ is said to be a $k$-independent set of $G$ if $\Delta(G[S])\leq k$, where $\Delta$ stands for the maximum degree. By the $k$-independence number of $G$, denoted by $\alpha_k(G)$, we mean the maximum cardinality of any $k$-independent set in $G$. An independent set is in fact a 0-independent set. Let $\alpha(G)$ stand for the maximum size of an independent set. We have $\alpha_0(G)=\alpha(G)$.

\noindent {\bf The outline of the paper is as follows.} In Section 2 we obtain some bounds for the independence number (Theorems \ref{boundzeta}, \ref{boundzeta-strong}, \ref{boundzeta-strong-2} and Proposition \ref{turan}). These bounds are achieved by a new greedy algorithm which we call Cheap-Greedy algorithm (Remark \ref{remark1} and Corollary \ref{coro1}). We show that equality in Theorem \ref{boundzeta} happens for a large class of graphs (Proposition \ref{equality}). In Section 3 we obtain a lower bound for $\alpha_1(G)$ (Theorem \ref{bound-1-indep}). This bound is achieved by a greedy type algorithm called 1-Cheap-Greedy algorithm (Remark \ref{remark3}). We show that our bound is better than the best known bounds for some families of graphs. In Section 4 we obtain a lower bound for $\alpha_2(G)$ (Theorem \ref{bound-2-indep}). This bound is achieved by a polynomial time algorithm (Remark \ref{remark4}). At the end of the paper we obtain a result for $\alpha_k(G)$ when $G$ is forest and $k$ is an arbitrary integer. Then we suggest some topics for further researches.

\section{Independence number}

\noindent The literature is full of papers concerning independence number in general and lower bounds for independence number in terms of vertex degrees, in particular. The history begins from the bounds of Caro \cite{C} and independently Wei \cite{W}. They proved that $\alpha(G)\geq {\sum}_{v\in G} 1/(deg_G(v)+1)$. The other types of lower bounds for independence number in terms of vertex degrees, generalizations or extensions of Caro-Wei bound can be found in \cite{ACL, G, H, HSch1, HSch2, M, STY, Se}. For some recent researches on this subject see \cite{ACL,H}. The main idea behind the Caro-Wei bound and also some other bounds is the following algorithm, the so-called Min-Greedy algorithm, which outputs an independent set of cardinality at least ${\sum}_{v\in G} 1/(deg_G(v)+1)$. The Min-Greedy algorithm repeatedly chooses a vertex say $v$ of the minimum degree in the underlying graph, puts it in the independent set and removes $N[v]$ from the graph. For convenience, denote ${\sum}_{v\in G} 1/(\zeta_G(v)+1)$ by $Z_1(G)$. The following bound improves the Caro-Wei bound and shows that the Min-Greedy algorithm outputs at least $Z_1(G)$ independent vertices.
It also shows that Min-Greedy algorithm not only can choose any vertex of minimum degree, at each step of the algorithm, but also it can choose any cheap vertex.

\begin{thm}
$$\alpha(G) \geq \sum_{v\in G} \frac{1}{\zeta(v)+1}.$$\label{boundzeta}
\end{thm}

\noindent \begin{proof}
We prove the theorem by induction on the order of graph. Recall that $Z_1(G) = {\sum}_{v\in G} 1/(\zeta(v)+1)$. Let $G=K_1$. Then $\alpha(G)=1=Z_1(G)$. Let now $G$ be an arbitrary graph and assume that the assertion holds for all graphs with less that $|G|$ vertices. As we observed earlier, there exists a cheap vertex $u$ in $G$ satisfying $\zeta(u) = deg_G(u)$ and $\zeta(u)= {\min}_{v\in N[u]}\zeta(v)$. Set $H=G\setminus N[u]$. By the induction hypothesis $\alpha(H) \geq Z_1(H)$. The following inequalities complete the proof.
$$Z_1(G) = \sum_{v\in V(H)} \frac{1}{\zeta(v)+1} + \sum_{v\in N[u]} \frac{1}{\zeta(v)+1}\leq Z_1(H)+  \sum_{v\in N[u]} \frac{1}{\zeta(u)+1}$$
$$\leq \alpha(H) + \frac{deg(u)+1}{deg(u)+1} = \alpha(H)+1 \leq \alpha(G).$$
\end{proof}

\noindent Based on the proof of Theorem \ref{boundzeta} and using the concept of cheap vertices, we define the {\it Cheap-Greedy algorithm}. We present the exact form of this algorithm after Theorem \ref{boundzeta-strong-2}, but its general format is as follows. In the following algorithm, by $CHEAP(G)$ we mean the set of all cheap vertices in $G$. Note that by Proposition \ref{output-cheap}, there exists an efficient algorithm which outputs $CHEAP(G)$ for any input graph $G$. We use this fact in line (3) of the following algorithm.

\noindent {\bf Cheap-Greedy algorithm for independence number (general format):}

\noindent (1) Input: graph $G$.

\noindent (2) Set $I=\varnothing$.

\noindent (3) Let $CHEAP(G)$ be the subgraph of $G$ induced by the cheap vertices of $G$.

\noindent (4) Pick any independent subset $S$ in $CHEAP(G)$ using Min-Greedy algorithm.

\noindent (5) Replace $I$ by $I\cup S$ and $G$ by $G\setminus N[S]$. If $V(G)\neq \varnothing$, then go to the step (3). Otherwise go to the step (6).

\noindent (6) Output the set $I$.

\noindent For example if $S$ is a single vertex in the above algorithm then according to the proof of Theorem \ref{boundzeta} we obtain an independent set of cardinality $Z_1(G)$. But since in general the set of cheap vertices in a graph is much more larger than the set of vertices with minimum degree then we should choose the set $S$ in the algorithm by a more intelligent method in order to obtain independent sets with more than $Z_1(G)$ vertices. Theorems \ref{boundzeta-strong} and \ref{boundzeta-strong-2} tell us how to choose such ``good" sets $S$ in each step of the execution of Cheap-Greedy algorithm. Note also that since any vertex with minimum degree is cheap then Min-Greedy is a particular case of Cheap-Greedy. But Example \ref{example} shows that Cheap-Greedy algorithm is not reduced to Min-Greedy algorithm. We have the following remark.

\begin{remark}
The Cheap-Greedy algorithm returns an independent set containing $Z_1(G)$ vertices. Moreover, Cheap-Greedy algorithm is not reduced to Min-Greedy algorithm.\label{remark1}
\end{remark}

\noindent The following example introduces a family of graphs for which Cheap-Greedy algorithm outputs a maximum independent set but any execution of Min-Greedy algorithm obtains a non-maximum independent set.

\begin{ex}\label{example}
Consider the graph $G_k$ whose vertex set is partitioned into $V(G)=I_1 \cup I_2 \cup \ldots \cup I_{2k}$ with $|I_i|=i$. For each $i$, $I_i$ is an independent set in $G$. For any $i$ with $1 \leq i \leq 2k-1$, the subgraph of $G$ induced on $I_i\cup I_{i+1}$ is complete bipartite. Let $I_1=\{u\}$. Then $u$ is the only vertex with minimum degree two in $G$. We observe that for any $i\in \{1, \ldots, 2k-2\}$ any vertex of $I_i$ has degenerate degree $i+1$. Any vertex in $I_{2k-1}\cup I_{2k}$ has degenerate degree $2k-1$. The set of cheap vertices in $G_k$ is $\{u\} \cup I_{2k}$, where the only vertex with minimum degree is $u$. It is seen that any execution of Min-Greedy algorithm outputs an independent set with no more than $1+3+5+\cdots+(2k-3)+(2k)$ vertices. But if we choose any cheap vertex in $I_{2k}$ (in fact Cheap-Greedy algorithm chooses the whole $I_{2k}$ at the first step of the algorithm) then the algorithm proceeds toward a maximum independent set of cardinality ${\sum}_{i=1}^k 2i$.\label{example}
\end{ex}

\noindent The following proposition shows that if at least one component of a graph $G$ is not regular, then the bound of Theorem \ref{boundzeta} is strictly better than the bound of Caro and Wei.

\begin{prop}
We have $Z_1(G)\geq {\sum}_{v\in G} \left( 1/[deg_G(v)+1] \right)$. And the equality holds if and only if each connected component of $G$ is regular.
\end{prop}

\noindent \begin{proof}
Without loss of generality we may assume that $G$ is a connected graph. If $G$ is regular of degree say $r$, then for any vertex $v$ of $G$, $\zeta(v)=deg(v)=r$. In this case both bounds have the same value. Assume that $Z_1(G) = {\sum}_{v\in G} 1/[deg_G(v)+1]$. Then $\zeta(v)=deg(v)$, for any vertex $v$. Let $j=\max \{\zeta(v): v\in G\}$. Let $u$ be any vertex of $G$ with $\zeta_G(u)=j$. By the definition of $\zeta(u)$, there exists a subgraph $H$ of $G$ such that $\delta(H)=j$ and $u\in H$. It is clear that $\zeta_G(v)=j$ for any vertex $v$ of $H$. Let $v$ be any vertex of $H$. We have $j=\delta(H)\leq deg_H(v) \leq deg_G(v) = \zeta_G(v) = j$. It follows that $H$ is a $j$-regular graph. If $G\neq H$ then there exist two adjacent vertices $x$ and $y$ such that $x\in H$ and $y\in G\setminus H$. Then $deg_G(x)>deg_H(x)=j$. We have now $j= \zeta_G(x)=deg_G(x)>j$. This contradiction shows that $G$ is regular.
\end{proof}

\noindent In Proposition \ref{turan} we obtain another Tur\'{a}n type bound for $\alpha(G)$. For any graph $G$ define the average degree of $G$ by ${\sum}_{v\in G} deg_G(v) /|G|$ and denote it by $\bar{d}$. A result of Tur\'{a}n \cite{Tu} asserts that if $G$ is a graph on $n$ vertices and with average degree $\bar{d}$, then $\alpha(G)\geq n/(\bar{d}+1)$.
In the following proposition we strengthen this result. In a graph $G$ set $\bar{\zeta}=({\sum}_{v\in G} \zeta(v))/n$. The following result is  obtained by Theorem \ref{boundzeta}.

\begin{prop}
$$\alpha(G)\geq \frac{n}{\bar{\zeta}+1}.$$\label{turan}
\end{prop}

\noindent \begin{proof}
By Theorem \ref{boundzeta} and the Cauchy-Schwartz inequality we have
$$\alpha(G) \geq \sum_{v\in G} \frac{1}{\zeta(v)+1} \geq \frac{n^2}{\sum_{v\in G}(\zeta(v)+1)}= \frac{n}{\bar{\zeta}+1}.$$
\end{proof}

\noindent Let $G$ be a graph with the maximum degree $\Delta(G)$. The Caro-Wei bound yields $\alpha(G)\geq n/(\Delta(G)+1)$. In this paragraph using Theorem \ref{boundzeta} we easily prove $\alpha(G)\geq (n/\Delta(G))-(2/3)$, provided that $G$ is not complete graph and odd cycle. First, let $G$ be a connected non-regular graph. Note that $\zeta(v)\leq \Delta(G)-1$, for any vertex $v$ of $G$, since otherwise $G$ should be regular. Theorem \ref{boundzeta} now implies that $\alpha(G)\geq n/\Delta(G)$. Let now $G$ be regular. Let $H$ be a graph obtained by $G$ by subdividing an edge of $G$. Then $\zeta(v)=\Delta(G)-1$ for any vertex of $G$ and $\zeta(u)=2$, where $u$ is the new vertex added to $G$. Applying Theorem \ref{boundzeta} for $H$ yields $\alpha(G) \geq (n/\Delta(G)) - (2/3)$.

\noindent As proved by Wei \cite{W}, the equality in Caro-Wei bound holds if and only if each connected component of $G$ is a complete graph. But in the bound of Theorem \ref{boundzeta}, equality holds for a larger family of graphs. In the following we introduce this family which we denote by ${\mathcal{F}}$. A typical member $G$ of ${\mathcal{F}}$ is obtained as follows. The vertex set of $G$ is decomposed into vertex disjoint complete subgraphs. In other words, $V(G)=C_1 \cup \ldots \cup C_k$, where the subgraph of $G$ induced by $C_i$ is a complete graph on say $n_i$ vertices. At this moment the degenerate degree of any vertex of $C_i$ is $n_i-1$, for each $i$. Now add arbitrary number of edges between these complete subgraphs of $G$ provided that the degenerate degree of no vertex is increased. We have now the following proposition.

\begin{prop}
Let $G$ be a graph. Equality holds for $G$ in the bound of Theorem \ref{boundzeta} if and only if $G\in {\mathcal{F}}$.\label{equality}
\end{prop}

\noindent \begin{proof}
Let $G\in {\mathcal{F}}$ and $V(G)=C_1 \cup \ldots \cup C_k$, where the subgraph of $G$ induced by $C_i$ is a complete graph on say $n_i$ vertices.
By the additional property of $G$ we have $\zeta(v)=n_i-1$, for any $v\in C_i$. We obtain $Z_1(G)=k$. From the other hand, $\alpha(G)\leq k$ and by Theorem \ref{boundzeta} $\alpha(G)\geq k$. This proves the one side of the proposition. Assume now that $Z_1(G)=\alpha(G)$. We prove the assertion by induction on the number of vertices. Let $u$ be a vertex of the minimum degree in $G$ and set $H=G\setminus N[u]$. Since $u$ is cheap in $G$, we have $Z_1(G)\leq 1+Z_1(H)\leq 1+\alpha(H)\leq \alpha(G)=Z_1(G)$. It follows that $\alpha(H)=Z_1(H)$. By applying the induction hypothesis for $H$ we obtain $H\in {\mathcal{F}}$. Hence the vertex set of $H$ can be partitioned into complete subgraphs as $V(H)= C_2 \cup \ldots \cup C_k$. Also for any $v\in H$ we have $\zeta_H(v)=\zeta_G(v)$.
Let $w\in H$ be a vertex such that $deg_H(w)=\delta(H)$. Assume that $w\in C_j$, for some $j\in \{2, \ldots, k\}$. We claim that $deg_{C_j}(w)=deg_H(w)$.
Otherwise $deg_{C_j}(w) < deg_H(w)$. It follows that $1+ \delta(H[C_j]) \leq deg_{C_j}(w) +1 \leq \delta(H)$. This contradicts $H\in {\mathcal{F}}$. Therefore $deg_{C_j}(w)=deg_H(w)$. It implies that $w$ is cheap in $H$. And since $\zeta_H(w)=\zeta_G(w)$, then $w$ is cheap in $G$. Hence we have $deg_H(w)=deg_G(w)=deg_{C_j}(w)$. Now let $H'=G\setminus N[w]$. Again we have $\alpha(H')=Z_1(H')$, since $w$ is cheap in $G$. Note that $G[C_j]$ is a complete subgraph and $deg_G(w)=deg_{C_j}(w)$. This mean that $G[N[w]]$ is a complete subgraph of $G$. Since $\alpha(H')=Z_1(H')$, we may repeat the same argument for $H'$ and deduce that $H'\in {\mathcal{F}}$. Consider the corresponding decomposition of $V(H')$ into complete subgraphs $D_1, \ldots, D_p$. These complete subgraphs of $H'=G\setminus N[w]$ together with $G[N[w]]$, introduce a decomposition of whole $G$ into complete subgraphs. Recall that $N[w]\subseteq C_j \subseteq H$ and $\zeta_{C_j}(v)=\zeta_H(v)=\zeta_G(v)$ for any $v\in N[v]$. Now $D_1, \ldots, D_p, N[w]$ have the required properties to conclude that $G\in {\mathcal{F}}$.
\end{proof}

\noindent The proof of Proposition \ref{equality} shows the following.

\begin{remark}
Let $G$ be any given graph. There exists a recursive algorithm that determines whether $G\in {\mathcal{F}}$.
\end{remark}

\noindent In the following result we strengthen Theorem \ref{boundzeta}. Let $G$ be a graph and $C$ be the set of cheap vertices of $G$. We need to define some subsets in $G$ and parameters associated with them.

\begin{itemize}
\item{Let $S$ be any maximal independent subset of $G[C]$. For example we can use any greedy algorithm based on the minimum degree to obtain $S$.}
\item{Let $B$ be the bipartite graph consisting of the partite sets $S$ and $N(S)$ and all edges of $G$ between $S$ and $N(S)$.}
\item{Let the connected components of $B$ be $B_1, B_2, \ldots, B_p$.}
\item{Denote the number of edges in $B_i$ by $e_i$.}
\item{Define $s_i=|B_i\cap S|$, $t_i=|B_i\cap N(S)|$ and $\lambda_i=1-\frac{e_i}{s_i}+\frac{t_i}{s_i}$.}
\end{itemize}

\noindent Note that for any $i$, $\lambda_i\leq 1$ and that for any distinct $i$ and $j$, $N[S\cap B_i]\cap N[S\cap B_j]=\varnothing$. We present the following theorem.

\begin{thm}
Let $G$ be any graph. Let $S$, $B_i$, $p$ and $\lambda_i$ be as defined above such that for any $i$, $\lambda_i\geq 0$. Then
$$\alpha(G) \geq \sum_{i=1}^p \sum_{v\in N[S\cap B_i]} \frac{1}{\zeta(v)+\lambda_i} + \sum_{v\in G\setminus N[S]} \frac{1}{\zeta(v)+1}.$$\label{boundzeta-strong}
\end{thm}

\noindent \begin{proof}
\noindent First note that $\lambda_i\geq 0$ means that $B_i$ is a connected graph with at most one cycle. Let $i$ be any arbitrary and fixed value with $1\leq i\leq p$. We show that
$$Z(B_i) := {\sum}_{v\in V(N[B_i\cap S])} \frac{1}{(\zeta(v)+\lambda_i)} \leq |B_i\cap S|.~~~~~\clubsuit$$
\noindent Let the vertices of $S\cap B_i$ be ordered such that $\zeta(u_1)\geq \ldots \geq \zeta(u_{s_i})$. The contribution of $u_1$ and its neighbors in
$Z(B_i)$ is at most $(\zeta(u_1)+1)/(\zeta(u_1)+\lambda_i)$. But the other vertices may have common neighbor with $u_1$. Since $B_i$ is connected we go from the neighbors of $N(u_1)$ to some other vertices in $B_i\cap S$ and determine their contribution in the above sum. Each typical contribution corresponding to such a vertex $u_j$ in $B_i\cap S$ is of the form $(\zeta(u_j)+1-r_j)/(\zeta(u_j)+\lambda_i)$, for some adequate $r_j\geq 1$ which is obtained by our method of determining the contribution of $N[u_j]$ in $Z(B_i)$. In fact $r_j$ is the number of the neighbors of $u_j$ which have been counted already as the neighbors of previous (with respect to $u_j$) vertices in $B_i\cap S$. For any constant $r\geq 0$, the real function $f(x)= (x-r)/(x+\lambda_i)$ is increasing in the interval $(0,\infty)$. Hence each typical element $(\zeta(u_j)+1-r_j)/(\zeta(u_j)+\lambda_i)$
is at most $(\zeta(u_1)+1-r_j)/(\zeta(u_1)+\lambda_i)$. We obtain the following inequality.
$$Z(B_i) \leq \frac{(\zeta(u_1)+1)+(\zeta(u_1)+1-r_2)+(\zeta(u_1)+1-r_3)+ \cdots + (\zeta(u_1)+1-r_{s_i})}{\zeta(u_1)+\lambda_i}.$$
\noindent We make the following claim.

\noindent {\bf Claim.} ${\sum}_j r_j=e_i-t_i$

\noindent The value in the left side of above equality is the number of repeated edges in counting the neighbors of $u_1, \ldots, u_{s_i}$. Let the vertices in $B_i\cap N(S)$ be $w_1, \ldots, w_{t_i}$. ``Each vertex $w_j$, has been counted $deg_{B_i}(w_j)-1$ many times repeatedly." Hence ${\sum}_j r_j= {\sum}_j (deg_{B_i}(w_j)-1)=e_i-t_i$. This proves the claim.

\noindent It follows that $$Z(B_i) \leq \frac{s_i\zeta(u_1)+s_i-{\sum}_j r_j}{\zeta(u_1)+\lambda_i}\leq \frac{s_i\zeta(u_1)+s_i-e_i+t_i}{\zeta(u_1)+\lambda_i}=s_i.$$

\noindent Set for simplicity $Z(S)={\sum}_{i=1}^p Z(B_i)$ and $Z=Z(S) + {\sum}_{v\in G\setminus N[S]} 1/(\zeta(v)+1).$
\noindent For each $i=1, \ldots, p$, the inequality $\clubsuit$ holds. Let $H=G\setminus N[S]$. Applying Theorem \ref{boundzeta} for $H$, we obtain $$Z\leq \sum_{i=1}^p s_i + \sum_{v\in H} \frac{1}{\zeta(v)+1} \leq |S| + \alpha(H) \leq \alpha(G).$$
\end{proof}

\noindent Denote the set of edges between $S$ and $N(S)$ in any graph by $E(S,N(S))$. Also by the edge density of a graph $H$ we mean $|E(H)|/|V(G)|$. Some graphs contain independent sets $S$ of cheap vertices such that the edge density in bipartite graph $S\cup N(S)$ is high (for example when $S$ and $N(S)$ form a complete bipartite graph with at least two cycles). In such situations the corresponding $\lambda$ i.e. $1+(|N(S)|-|E(S,N(S))|)/|S|$ is negative and hence we cannot use Theorem \ref{boundzeta-strong} (instead we apply Theorem \ref{boundzeta-strong-2}). In fact in order to have better lower bounds for $\alpha(G)$ we need smaller (and then possibly negative) $\lambda$. We need some explanations before stating the next theorem. Let $S$ be any independent set of cheap vertices with the same degenerate degree in a graph $G$. For any $S'\subseteq S$ define $\lambda(S')=1+[(|N(S')|-|E(S',N(S'))|)/|S'|]$. In order to optimum use of Theorem \ref{boundzeta-strong-2} and also run Cheap-Greedy algorithm we need to choose a subset $S_0\subseteq S$ such that $\lambda(S_0)$ is the smallest possible value among the other subsets $S'$ of $S$. Let $B(S)$ be the bipartite graph on the partite sets $S$ and $N(S)$ and all edges of $G$ between $S$ and $N(S)$. One way to fulfil this purpose is to choose $S_0 \subseteq S$ such that $S_0 \cup N(S_0)$ has the maximum edge density in $B(S)$ (call this method the maximum edge density method). In \cite{P}, the author presents a polynomial time algorithm such that given a graph $G$, it returns a subgraph of $G$ with the maximum edge density. Note also that in any graph $G$ we can obtain all vertices with the same degeneracy degree say $i$ by the algorithm presented in Proposition \ref{algor}.

\begin{thm}
Let $G$ be any graph with the degeneracy $p$. For any $i\in \{\delta(G), \ldots, p\}$, let $C_i$ be the set of cheap vertices with degenerate degree $i$. Let $S_i$ be any independent set in $G[C_i]$ (preferably by the method of maximum edge density). Define $\lambda_i= 1-(e_i/|S_i|)+ (|N(S_i)|/|S_i|)$, where $e_i$ is the number of edges between $S_i$ and $N(S_i)$. Let $j$ be such that $\lambda_j=\min_i \lambda_i$. Set $\lambda=\lambda_j$, $S=S_j$ and $e=e_j$. Then $$\alpha(G) \geq \sum_{v\in N[S]} \frac{1}{\zeta(v)+\lambda} + \sum_{v\in G\setminus N[S]} \frac{1}{\zeta(v)+1}.$$\label{boundzeta-strong-2}
\end{thm}

\noindent \begin{proof}
Similar to the proof of Theorem \ref{boundzeta-strong}, we show that
$$Z(S) := {\sum}_{v\in V(N[S])} \frac{1}{(\zeta(v)+\lambda)} \leq |S|.$$
\noindent Let $S=\{u_1, u_2, \ldots, u_{|S|}\}$. The contribution of $u_1$ and its neighbors in
$N(S)$ is at most $(\zeta(u_1)+1)/(\zeta(u_1)+\lambda)$. The contribution of $N[u_2]$ is at most $(\zeta(u_2)+1-r_2)/(\zeta(u_2)+\lambda)$, where $r_2$ is the number of common neighbors of $u_1$ and $u_2$ in $N(S)$. We continue to determine the contribution of the other vertices of $S$. For any vertex $u_p$ the contribution is of the form $(\zeta(u_p)+1-r_p)/(\zeta(u_p)+\lambda)$, for some $r_p\geq 1$. We obtain the following
$$Z(S) \leq \frac{(\zeta(u_1)+1)+(\zeta(u_1)+1-r_2)+(\zeta(u_1)+1-r_3)+ \cdots + (\zeta(u_1)+1-r_{b_i})}{\zeta(u_1)+\lambda}.$$
\noindent We have the equality ${\sum}_p r_p=e-|N(S)|$. Its proof is similar to the proof of the claim in Theorem \ref{boundzeta-strong}. It follows that $$Z(S) \leq \frac{|S|\zeta(u_1)+|S|-{\sum}_p r_p}{\zeta(u_1)+\lambda} = \frac{|S|\zeta(u_1)+|S|-e+|N(S)|}{\zeta(u_1)+\lambda}=|S|.$$
\noindent The rest of the proof is similar to the proof of Theorem \ref{boundzeta-strong}. We omit the details.
\end{proof}

\noindent We end this section by presenting the complete form of Cheap-Greedy algorithm. Recall the notations $S$, $Z(S)$ and $\lambda_1, \ldots, \lambda_p$ from the proof of Theorem \ref{boundzeta-strong}. Define $\lambda = \min_i \lambda_i$. Also recall $S$, $\lambda$ and $Z(S)$ from the proof of Theorem \ref{boundzeta-strong-2}. We use these sets and values in the following algorithm.

\noindent {\bf Cheap-Greedy algorithm for independence number:}

\noindent (1) Input: graph $G$.

\noindent (2) Set $I=\varnothing$ and $Z=0$.

\noindent (3) Let $S_1$, $\lambda(S_1)$ and $Z(S_1)$ (resp. $S_2$, $\lambda(S_2)$ and $Z(S_2))$ be the corresponding values $S$, $\lambda$ and $Z(S)$ as in the proof of Theorem \ref{boundzeta-strong} (resp. Theorem \ref{boundzeta-strong-2}).

\noindent (4) If $\lambda(S_2)< \lambda(S_1)$ then replace $I$ by $I\cup S_2$, $Z$ by $Z+Z(S_2)$ and $G$ by $G\setminus N[S_2]$.

\noindent (5) If $\lambda(S_1) \leq \lambda(S_2)$ then replace $I$ by $I\cup S_1$, $Z$ by $Z+Z(S_1)$ and $G$ by $G\setminus N[S_1]$.

\noindent (6) Output $I$ and $Z$.

\noindent We have the following corollary.

\begin{cor}
Cheap-Greedy algorithm outputs an independent set with at least $Z$ vertices.\label{coro1}
\end{cor}

\section{1-independent sets}

\noindent In this section we study 1-independent sets. The $k$-independent sets and related topics were widely studied in the literature, see e.g. \cite{CH, CT, F, HS}. Most of the studies in the area of $k$-independent sets focus on lower bounds for $\alpha_k(G)$. The best known bound for $\alpha_1(G)$ in terms of average degree is from \cite{CH}, $\alpha_1(G)\geq 2n/(\lceil \bar{d}(G) \rceil+2)$, where $\bar{d}(G)=({\sum}_{v\in G} deg_G(v))/|G|$. Also Caro and Tuza \cite{CT} obtained the bound $\alpha_1(G)\geq {\sum}_{v\in G} 3/[2(deg(v)+1)]$ in terms of the degree sequence of $G$ with $\delta(G)\geq 2$. Generalize the notation $Z_1(G)$ and define $Z_k(G)$ as follows.
$$Z_k(G):= \sum_{v\in G} \min \left\{1, \frac{1}{\zeta(v)+(1/k)}\right\}.$$
\noindent Let $V_0$ be the set of isolated vertices of $G$. We observe that
$$Z_k(G)= \left(\sum_{v\in V(G)\setminus V_0} \frac{1}{\zeta(v)+(1/k)}\right) +|V_0|.$$ Note that when $k=1$ then $({\sum}_{v\in V(G)\setminus V_0} 1/(\zeta(v)+1))+|V_0| = {\sum}_{v\in V(G)} 1/(\zeta(v)+1)=Z_1(G)$. In the following we prove that $\alpha_1(G)\geq Z_{2}(G)$. We first introduce the concept of $k$-cheap sets. Then we use $1$-cheap sets in proving our bound for $\alpha_1(G)$. By a $k$-cheap set in a graph $G$ we mean either any subset of isolated vertices in $G$ or any subset $S$ of non-isolated vertices of $G$ such that the following conditions hold.

{\bf (i)} $\sum_{v\in N[S]} \frac{1}{\zeta(v)+\frac{1}{k+1}} \leq |S|$,~~~~~~~~~~ {\bf (ii)} $\Delta(G[S]) \leq k.$

\noindent The reason why we call a set $S$ satisfying the above conditions a $k$-cheap set, is that the contribution of $N[S]$ in $Z_{k+1}(G)$ is at most $|S|$, in other words $S$ has low cost in $Z_{k+1}(G)$. For each $k$ the minimal $k$-cheap subsets have important role in proving our bounds for $\alpha_k(G)$. For $k=0$ the only minimal cheap subset is a single vertex. This is equivalent to our earlier definition of cheap vertex. The conditions (i) and (ii) are easily verified for any single cheap vertices.

\noindent First, we investigate the structure of 1-cheap sets. There are three minimal isolated-vertex-free 1-cheap sets which we call type I, type II and type III. The type I consists of two adjacent cheap vertices. The type II consists of two adjacent vertices $u$ and $w$ such that $u$ is cheap in $G$, $u$ is the only cheap vertex in $G$ which is adjacent to $w$; and $w$ is a cheap vertex in $G\setminus u$. The type III consists of two non-adjacent cheap vertices such that they have at least one common neighbor. These 1-cheap sets are minimal in the sense that any 1-cheap set in a graph contains at least one of them (Theorem \ref{1-cheap}). We show that a subset $S$ of type II is a 1-cheap set. The proofs for two other types are similar. Let $V(S)=\{u,v\}$. It is enough to verify the condition (i) in the definition of cheap sets. Note that $\zeta(u)=deg(u)$, $\zeta(w)=deg(w)-1$ and $\zeta(u)\leq \zeta(w)$. We have $$\sum_{v\in N[S]} \frac{1}{\zeta(v)+\frac{1}{2}} \leq \frac{deg(u)}{deg(u)+(1/2)} + \frac{\zeta(w)+1}{\zeta(w)+(1/2)} \leq \frac{2deg(u)+1}{deg(u)+(1/2)} =2.$$
\noindent The following theorem shows that the 1-cheap sets of type I, II and III are in fact the only minimal 1-cheap sets, because any graph $G$ contains at least one of them.

\begin{thm}
Any graph $G$ with at least one edge contains a $1$-cheap set of type I or II and or III. Moreover, such a set can be found by a greedy algorithm.\label{1-cheap}
\end{thm}

\noindent \begin{proof}
Let $I$ be the set of non-isolated cheap vertices in $G$. If $I$ is not an independent set then we obtain two cheap vertices which are adjacent. These two vertices form a 1-cheap set of type I. If there are two non-adjacent cheap vertices in $I$ with a common neighbor in $V(G)\setminus I$, then these two cheap vertices form a 1-cheap set of type III. Assume hereafter that $I$ is an independent set and no two elements of $I$ have a common neighbor. Set $H=G\setminus I$. Let $v$ be any cheap vertex in $H$. We claim that $v$ has a neighbor in $I$. Otherwise, $deg_G(v)=deg_H(v)$. It is clear that $\zeta_H(v)\leq \zeta_G(v)$. We have also $deg_H(v)=\zeta_H(v)$, since $v$ is cheap in $H$. This shows that $deg_G(v)=deg_H(v) =\zeta_H(v) \leq \zeta_G(v)$. It follows that $deg_G(v) \leq \zeta_G(v)$. We have also $deg_G(v) \geq \zeta_G(v)$ by the definition of $\zeta_G(v)$. Hence $deg_G(v) = \zeta_G(v)$. Note also that for any vertex $u\in H$, $\zeta_H(u)\leq \zeta_G(u)$. Let $N[v]$ be the set of neighbors of $v$ in $G$ together with the vertex $v$ itself. Since $v$ is cheap in $H$, then $$\zeta_G(v)=\zeta_H(v)={\min}_{u\in N[v]\cap V(H)}\zeta_H(u)={\min}_{u\in N[v]}\zeta_H(u)\leq {\min}_{u\in N[v]}\zeta_G(u).$$
\noindent It follows that $v$ is cheap in $G$. This contradiction shows that $v$ is adjacent to only one cheap vertex say $w$ in $G$. Now, two vertices $v$ and $w$ form a 1-cheap set of type II.

\noindent Let us now present an algorithm. We first search for a pair $(u,w)$ of cheap vertices such that either $uw\in E(G)$ or $u$ and $w$ have a common neighbor. If we find no such pair then we remove all cheap vertices from $G$ and obtain a cheap vertex in the remaining graph. This cheap vertex is adjacent to a vertex which is cheap in $G$. By this procedure we obtain a 1-cheap set.
\noindent \end{proof}

\noindent Recall that $Z_2(G) = ({\sum}_{v\in V(G)\setminus V_0} 1/(\zeta(v)+(1/2))) +|V_0|$, where $V_0$ is the set of isolated vertices in $G$. We have now the following theorem.

\begin{thm}
$$\alpha_1(G) \geq Z_2(G).$$\label{bound-1-indep}
\end{thm}

\noindent \begin{proof}
We prove the theorem by induction on the order of graph. It is valid for $K_1$. Let $G$ be an arbitrary graph and assume that the assertion holds for all graphs with order less that $|G|$. Let $V_0$ be the set of isolated vertices in $G$. By Theorem \ref{1-cheap} there exists a $1$-cheap set say $S$ in $G\setminus V_0$. Set $H=G\setminus (N[S]\cup V_0)$. Note that for any $v\in H$, $\zeta_H(v) \leq \zeta_G(v)$ and for any $v\in N[S]$, $1/(\zeta(v)+(1/2)) < 1$. We apply the induction assertion for $H$, use the fact that $S$ is $1$-cheap in $G\setminus V_0$ and obtain the following inequalities which complete the proof.
$$Z_{2}(G) = |V_0| + \sum_{v\in V(H)} \min \{1, \frac{1}{\zeta(v)+(1/2)}\}+ \sum_{v\in N[S]} \min \{1, \frac{1}{\zeta(v)+(1/2)}\}$$
$$\hspace{-4cm}\leq |V_0| + Z_{2}(H) + \sum_{v\in N[S]} \frac{1}{\zeta(v)+(1/2)}$$
$$\hspace{-6cm}\leq |V_0| + \alpha_1(H) + |S|$$
$$\hspace{-7.8cm}\leq \alpha_1(G).$$
\noindent \end{proof}

\noindent Based on Theorems \ref{1-cheap} and \ref{bound-1-indep} we present the following greedy algorithm for $\alpha_1(G)$. The algorithm is greedy type because in each step of the execution of the algorithm, it searches for a set with low cost i.e. 1-cheap set.\\

\noindent {\bf 1-Cheap-Greedy algorithm for 1-independence number:}

\noindent (1) Input: graph $G$.

\noindent (2) Set $I=\varnothing$.

\noindent (3) Let $V_0(G)$ be the set of isolated vertices in $G$. Pick a 1-cheap set $S$ in $G\setminus V_0(G)$ using the greedy algorithm presented in Theorem \ref{1-cheap}.

\noindent (4) Replace $I$ by $I \cup V_0(G) \cup V(S)$ and $G$ by $G\setminus (N[S]\cup V_0(G))$. If $V(G)\neq \varnothing$, then go to the step (3). Otherwise go to the step (5).

\noindent (5) Output $G[I]$.\\

\noindent We have the following remark by Theorem \ref{bound-1-indep}.

\begin{remark}
Let 1-cheap-Greedy algorithm output the set $I$ for the graph $G$. Then $G[I]$ is 1-independent set of cardinality at least $Z_2(G)$.\label{remark3}
\end{remark}

\noindent Note that the bound of Theorem \ref{bound-1-indep} is tight for paths. Let $G_k$ be the path on $3k$ vertices, where $k$ is any positive integer. It can be easily checked that $\alpha_1(G_k)=2k$. Hence, in the bound of Theorem \ref{bound-1-indep} equality holds for $G_k$. The best known bound for $\alpha_1(G)$ is $2n/(\bar{d}(G)+2)$ from \cite{CH}, where $n$ and $\bar{d}(G)$ are the order and average degree of $G$, respectively. We note that if $\bar{\zeta}(G) < (\bar{d}(G)+1)/2$ then the bound of Theorem \ref{bound-1-indep} is better than the bounds of \cite{CH} and \cite{CT}. There are many graphs satisfying the latter inequality. For example edge maximal graphs with a given degenerate degree sequence (including edge maximal $\zeta$-regular graphs) satisfy this condition.

\section{2-independent sets}

\noindent In this section we prove that $\alpha_2(G)\geq Z_3(G)$. The proof is based on the concept of 2-cheap sets. Let $C_1$ be the set of cheap vertices in $G$. Let $C_2$ be the set of cheap vertices in $G\setminus C_1$. Similarly, define $C_i$ as the set of cheap vertices in $G\setminus (C_1 \cup \ldots \cup C_{i-1})$. It is clear that the vertex set of $G$ is partitioned into $C_1, \ldots, C_t$ for some integer $t$. We begin with the following lemma. In Theorem \ref{2-cheap} we prove that any graph has a 2-cheap set. Our scenario to prove this theorem is based on contradiction. So in the following lemmas we assume (on the contrary) that an arbitrary graph $G$ contains no 2-cheap set and explore the properties of such graphs $G$ and then obtain the desired contradiction in Theorem \ref{2-cheap}. Note that since we have assumed in the following lemmas that $G$ has not any 2-cheap sets then $G$ does not contain any isolated vertices. Hence we have $Z_2(G)= {\sum}_{v\in G} 1/ (\zeta(v)+(1/3))$. This mean that in the following lemmas we take into account this simplified expression of $Z_2(G)$.

\begin{lemma}
Let $G$ be a graph without any 2-cheap set. Then

\noindent (i) $C_1$ is an independent set.

\noindent (ii) No three vertices of $C_1$ has a common neighbor. In particular no vertex in $C_2$ has three neighbors in $C_1$.

\noindent (iii) Any vertex of $C_2$ has exactly one neighbor in $C_1$.

\noindent (iv) For any $v\in C_2$ we have $\zeta(v)=deg(v)-1$.

\noindent (v) Any vertex of $C_1$ has at most one neighbor in $C_2$.

\noindent (vi) $C_2$ is an independent set.\label{yek}
\end{lemma}

\noindent \begin{proof}
To prove $(i)$, assume on the contrary that there exist two adjacent vertices $u$ and $w$ in $C_1$. Set $S=\{u,v\}$. Note that $deg(u)=\zeta(u)$ and $deg(w)=\zeta(w)$. We have
$$\sum_{v\in N[S]} \frac{1}{\zeta(v)+\frac{1}{3}} \leq \frac{deg(u)}{deg(u)+(1/3)} + \frac{deg(w)}{\zeta(w)+(1/3)} \leq 2.$$
\noindent It follows that $S$ is a 2-cheap set, a contradiction.

\noindent To prove $(ii)$, let some three vertices $u$, $w$ and $p$ from $C_1$ have a common neighbor. Set $S=\{u, w, p\}$. Assume that $\zeta(u)\geq \zeta(w) \geq \zeta(p)$. The following yields a contradiction.
$$\sum_{v\in N[S]} \frac{1}{\zeta(v)+\frac{1}{3}} \leq \frac{deg(u)+1}{deg(u)+(1/3)} + \frac{deg(w)}{\zeta(w)+(1/3)} + \frac{deg(p)}{\zeta(p)+(1/3)}\leq \frac{3\zeta(u)+1}{\zeta(u)+(1/3)}=3.$$

\noindent To prove $(iii)$, let $S$ be consisted of two non-adjacent cheap vertices $u$ and $w$ and their common neighbor $p$ in $C_2$. Note that by $(ii)$, $deg(p)\leq \zeta(p)+2$. We have
$$\sum_{v\in N[S]} \frac{1}{\zeta(v)+\frac{1}{3}} \leq \frac{deg(u)}{deg(u)+(1/3)} + \frac{\zeta(w)}{\zeta(w)+(1/3)} + \frac{deg(p)-1}{\zeta(p)+(1/3)}$$
$$\leq \frac{2deg(u)}{deg(u)+(1/3)} + \frac{\zeta(p)+1}{\zeta(p)+(1/3)} \leq \frac{3deg(u)+1}{deg(u)+(1/3)}=3.$$
\noindent It follows that $S$ is a 2-cheap set, a contradiction.

\noindent Part $(iv)$ is proved by the membership of $v$ in $C_2$ and also $(ii)$. To prove $(v)$, let $u\in C_1$ have two neighbors $w$ and $p$ in $C_2$. By $(iv)$, $\zeta(w)=deg(w)-1$ and $\zeta(p)=deg(p)-1$. We can easily show that $\{u, w, p\}$ forms a 2-cheap set. Finally we prove $(vi)$. Let $u, w\in C_2$ be two adjacent vertices. There exist $x, y\in C_1$ such that $ux, wy \in E(G)$. Consider the induced path on $S=\{x, u, w, y\}$. Similar to the previous cases it can be shown that $S$ is 2-cheap. We omit the details.
\end{proof}

\noindent We also need the following lemma.

\begin{lemma}
Let $G$ be a graph without a 2-cheap set. Let $P= u_1, u_2, \ldots u_m$ and $Q= w_1, w_2, \ldots, w_m$ be two vertex disjoint paths in $G$ such that for each $i=1, \ldots, m$, $\{u_i, w_i\} \subseteq C_i$. Suppose there exists no edge between these two paths and in particular $deg(v)=\zeta(v)+1$ for any $v\in V(P)\cup V(Q) \setminus \{u_1, w_1\}$. Then no vertex in $C_{m+1}$ is adjacent to $u_m, w_m$.\label{two-paths}
\end{lemma}

\noindent \begin{proof}
Let $\zeta(u_1)\geq \zeta(w_1)$. Assume on the contrary that $u_m$ and $w_m$ have a common neighbor. We prove that the set $S=V(P)\cup V(Q)$ is a 2-cheap set. In the following counting, we count each vertex $u_i$ (resp. $w_i$) as a vertex of $N[u_i]$ (resp. $N[w_i]$). It means that $u_i$ is counted neither in $N[u_{i-1}]$ nor in $N[u_{i+1}]$ and the similar holds for $w_i$. We obtain the following inequalities
$$\sum_{v\in N[S]} \frac{1}{\zeta(v)+\frac{1}{3}} \leq \frac{deg(u_1)}{\zeta(u_1)+(1/3)}+ \sum_{i=2}^{m-1} \frac{deg(u_i)-1}{\zeta(u_i)+(1/3)}
+ \frac{deg(u_m)}{\zeta(u_i)+(1/3)}$$
$$+ \frac{deg(w_1)}{\zeta(w_1)+(1/3)}+ \sum_{i=2}^{m} \frac{deg(w_i)-1}{\zeta(w_i)+(1/3)}$$
$$\leq \frac{2\zeta(u_1)}{\zeta(u_1)+(1/3)}+ \sum_{i=2}^{m-1} \frac{\zeta(u_i)}{\zeta(u_i)+(1/3)}+ \frac{\zeta(u_m)+1}{\zeta(u_m)+(1/3)}+ \sum_{i=2}^{m} \frac{\zeta(w_i)}{\zeta(w_i)+(1/3)}.$$
\noindent Since the real function $f(x)=(x+1)/(x+(1/3))$ is decreasing for $x\in (0, \infty)$, then $(\zeta(u_m)+1)/(\zeta(u_m)+(1/3))$ is at most $(\zeta(u_1)+1)/(\zeta(u_1)+(1/3))$. Hence $$\sum_{v\in N[S]} \frac{1}{\zeta(v)+\frac{1}{3}} \leq \frac{3\zeta(u_1)+1}{\zeta(u_1)+(1/3)} + (2m-3) = 2m.$$
\noindent It follows that $S=V(P)\cup V(Q)$ is a 2-cheap set. This contradiction proves the lemma.
\end{proof}

\noindent The following lemma provides more information on $C_3$.

\begin{lemma}
Let $G$ be a graph without a 2-cheap set. Then

\noindent (i) No vertex in $C_3$ is adjacent to more than one vertex in $C_2$.

\noindent (ii) There exists no edge between $C_1$ and $C_3$.\label{do}
\end{lemma}

\noindent \begin{proof}
Assume on the contrary that $(i)$ does not hold. Then we obtain two edges $u_1u_2$ and $w_1w_2$ such that $u_1, w_1\in C_1$, $u_2, w_2\in C_2$ and $u_2$ and $w_2$ have a common neighbor. Since $u_1$ and $w_1$ are cheap in $G$ then $\zeta(u_1)=deg_G(u_1)$ and $\zeta(w_1)=deg_G(w_1)$. Since $u_2, w_2\in C_2$ then by the part (iv) of Lemma \ref{yek} we have $\zeta(u_2)=deg_G(u_2)-1$ and $\zeta(w_2)=deg_G(w_2)-1$. Part (iv) of Lemma \ref{yek} shows that $S=\{u_1, u_2, w_1, w_2\}$ is an induced matching consisting of two edges with a common neighbor for $u_2$ and $w_2$ outside the set $S$. Lemma \ref{two-paths} implies that $S$ is a 2-cheap set. This contradiction proves part $(i)$.

\noindent We prove $(ii)$. Assume on the contrary that there exists $u\in C_1$ and $w\in C_3$ such that $u$ is adjacent to $w$. Since $w$ is not cheap in $G\setminus C_1$ then it has a neighbor say $y$ in $C_2$. Similarly, $y$ has a neighbor say $u'$ in $C_1$. There are two possibilities.

\noindent Case 1. $u=u'$. In order to obtain a contradiction we show that $T=\{u, y\}$ forms a 2-cheap set. In the following counting we count the vertex $y$ as an element of $N[y]$. This means that $y$ is not counted as an element in $N[u]$. Also we count $w$ only as a member of $N[u]$. Recall that $deg_G(y)-1=\zeta(y)$. The following inequalities show that $T$ is a 2-cheap set in $G$, i.e. a contradiction.
$$\sum_{v\in N[T]} \frac{1}{\zeta(v)+\frac{1}{3}} \leq \frac{deg(u)}{deg(u)+(1/3)} + \frac{deg(y)-1}{\zeta(y)+(1/3)} \leq 2.$$
\noindent Case 2. $u\not=u'$. In this case too, we show that $T'=\{u, u', y\}$ is a 2-cheap set. In the following counting, we count $w$ as an element of $N(y)$ and $y$ as an element of $N[y]$. This means that $w$ is not counted in $N[u]$ and $y$ is not counted in $N[u']$. Also $u'$ is not counted in $N[y]$. We have now the following inequalities.
$$\sum_{v\in N[T']} \frac{1}{\zeta(v)+\frac{1}{3}} \leq \frac{deg(u)}{deg(u)+(1/3)} + \frac{deg(u')}{deg(u')+(1/3)} + \frac{deg(y)}{\zeta(y)+(1/3)}$$
$$\leq \frac{deg(u)}{deg(u)+(1/3)} + \frac{deg(u)}{deg(u)+(1/3)} + \frac{\zeta(y)+1}{\zeta(y)+(1/3)}$$
$$\hspace{-4.8cm}\leq \frac{3deg(u)+1}{deg(u)+(1/3)}=3.$$
\end{proof}

\noindent The following lemma provides enough information on graphs having no 2-cheap sets. We need to define the concept of jumping edge, which is any edge consisting of two vertices $u$ and $w$ such that $u\in C_i$ and $w\in C_j$ for some $i$ and $j$ with $|i-j|\geq 2$.

\begin{lemma}
Let $G$ be a graph without a 2-cheap set. Let $i$ be any integer with $3\leq i\leq t$. Let also $H_i$ be the subgraph of $G$ induced by $C_1\cup \ldots \cup C_i$. Then $H_i$ does not contain any jumping edge.\label{jump}
\end{lemma}

\noindent \begin{proof}
We prove the lemma by the induction on $i$. Lemma \ref{do} proves the assertion for the case $i=3$. Assume that the assertion holds for $C_j$. We prove the assertion for $C_{j+1}$. Let $u$ be any vertex in $C_{j+1}$.

\noindent {\bf Claim 1.} The vertex $u$ is not adjacent to any two vertices in $C_j$.

\noindent {\bf Proof of Claim 1.} Assume on the contrary that Claim 1 is not valid. Then we obtain two paths $P= u_1, u_2, \ldots u_j$ and $Q= w_1, w_2, \ldots, w_j$ such that for each $i=1, \ldots, j$, $\{u_i, w_i\} \subseteq C_i$ and $u$ is adjacent to $u_j$ and $w_j$.
By the induction hypothesis there exists no edge between $V(P)$ and $V(Q)$. This implies that $deg(v)=\zeta(v)+1$ for any $v \in V(P)\cup V(Q) \setminus \{u_1, w_1\}$. Now by Lemma \ref{two-paths}, no vertex in $C_{j+1}$ is adjacent to $u_j$ and $w_j$. This contradiction proves Claim 1.

\noindent We make another claim for the arbitrary vertex $u$ of $C_{j+1}$.

\noindent {\bf Claim 2.} There exists no $r$ and $s$ such that $u$ is adjacent to $C_r$ and $C_s$.

\noindent {\bf Proof of Claim 2.} Let $r\leq s$. Assume on the contrary that the claim does not hold. Then we obtain two paths $P= u_1, u_2, \ldots u_r$ with $u_i\in C_i$ and $Q= w_1, w_2, \ldots, w_s$ with $w_i\in C_i$ such that $u_r$ and $w_s$ are adjacent to $u$. Again by the induction hypothesis there exists no edge between $V(P)$ and $V(Q)$. We have also $deg(v)=\zeta(v)+1$ for any $v \in V(P)\cup V(Q) \setminus \{u_1, w_1\}$. Now similar to the proof of Lemma \ref{two-paths} we obtain a contradiction. This completes the proof Claim 2.

\noindent The vertex $u$ has a neighbor say $u_j\in C_j$. There exists a path $P= u_1, u_2, \ldots u_j$ with $u_i\in C_i$. The only remaining case is to prove that for each $i$ with $1\leq i\leq j-1$, $u$ is not adjacent to $u_i$. Assume on the contrary that $u$ is adjacent to some $u_r$. This mean that $u_r$ and $u_j$ have a common neighbor. In this situation we show that the path $P$ forms a 2-cheap set. We use the same counting technique as in the proof of Lemma \ref{two-paths}. By this method the contribution of each $u_i$ ($i\leq j-1$) is at most $\zeta(u_i)/(\zeta(u_i)+(1/3))$. For the last vertex i.e. $u_j$ the contribution is potentially $(\zeta(u_j)+1)/(\zeta(u_j)+(1/3))$. But $u_j$ and $u_r$ have a common neighbor and this neighbor is already counted. It implies that the contribution of $u_j$ reduces to $\zeta(u_j)/(\zeta(u_j)+(1/3))$. This completes the proof.
\end{proof}

\noindent We are now ready to prove our existence theorem for 2-cheap sets.

\begin{thm}
Any graph $G$ contains a $2$-cheap set.\label{2-cheap}
\end{thm}

\noindent \begin{proof}
Let $G$ be a graph without any 2-cheap set. We set up the sets $C_1, \ldots, C_t$ as we explained before and partition the vertex set of $G$ into $C_1 \cup \ldots \cup C_t$. We prove by induction on $j$ that $C_j$ is an independent set. Lemma \ref{yek} proves this assertion for $j=1$. Assume that $C_1, C_2, \ldots, C_j$ are independent sets. We show that $C_{j+1}$ is independent. Let $u, w\in C_{j+1}$ but there is an edge between them. There exist two paths $P= u_1, u_2, \ldots u_j$ and $Q= w_1, w_2, \ldots, w_j$ with $\{u_i, w_i\} \subseteq C_i$, for each $i=1, \ldots, j$ such that $u$ (resp. $w$) is adjacent to $u_j$ (resp. $w_j$). We claim that a new path $R$ consisting of the vertices $u_1, u_2, \ldots u_j, u, w, w_1, w_2, \ldots, w_j$ forms a 2-cheap set. First note that Lemma \ref{jump} implies $deg(v)-1=\zeta(v)$ for any $v\in V(R)\setminus \{u_1, w_1\}$. The proof that $R$ is 2-cheap set is obtained similar to the previous methods. We conclude finally that $G$ is a disjoint vertex union of paths. It follows that $G$ itself is 2-cheap. This contradiction proves the theorem.
\end{proof}

\noindent The following theorem is easily proved.

\begin{thm}
$$\alpha_2(G) \geq Z_3(G).$$\label{bound-2-indep}
\end{thm}

\noindent \begin{proof}
We prove the theorem by induction on the order of graph. It is valid for $K_1$. Let $G$ be an arbitrary graph and assume that the assertion holds for all graphs with order less that $|G|$. Let $V_0$ be the set of isolated vertices in $G$. Let $V_0$ be the set of isolated vertices in $G$. By Theorem \ref{2-cheap} there exists a $2$-cheap set say $S$ in $G\setminus V_0$.
Set $H=G\setminus (N[S]\cup V_0)$. Note that for any $v\in H$, $\zeta_H(v) \leq \zeta_G(v)$ and for any $v\in N[S]$, $1/(\zeta(v)+(1/3)) < 1$.
We apply the induction assertion for $H$, use the fact that $S$ is $2$-cheap in $G\setminus V_0$ and obtain the following inequalities which complete the proof.
$$Z_{3}(G) = |V_0| + \sum_{v\in V(H)} \min \{1, \frac{1}{\zeta(v)+(1/3)}\}+ \sum_{v\in N[S]} \min \{1, \frac{1}{\zeta(v)+(1/3)}\}$$
$$\hspace{-4cm}\leq |V_0| + Z_{3}(H) + \sum_{v\in N[S]} \frac{1}{\zeta(v)+(1/3)}$$
$$\hspace{-6.7cm}\leq |V_0| + \alpha_2(H) + |S|$$
$$\hspace{-8.5cm}\leq \alpha_2(G).$$
\end{proof}

\noindent The bound of Theorem \ref{bound-2-indep} is better than the bound $\alpha_2(G)\geq 3n/(\bar{d}(G)+3)$ from \cite{CH} for some graphs, for example for edge maximal graphs with $\bar{\zeta}(G)<2$. This condition is equivalent to the following. Assume that there are $n_i$ vertices in $G$ with degenerate degree $i$ and let $1\leq i \leq k$. If ${\sum}_{i=3}^k (i-2)n_i < n_1$ then our bound for $\alpha_2(G)$ is strictly better then the bound $\alpha_2(G)\geq 3n/(\bar{d}(G)+3)$. From the other side, the proofs of lemmata and theorems in this section show that we can obtain a 2-cheap set by a polynomial time algorithm. These 2-cheap sets are either induced paths or cycles and or two induced paths whose endpoints have a common neighbor.

\begin{remark}
There exists a polynomial time algorithm which obtains a 2-independent set of cardinality at least $Z_2(G)$.\label{remark4}
\end{remark}

\section{Concluding remarks}

\noindent In this section we first give a result for $k$-independence number of forests when $k$ is an arbitrary positive integer and then we propose three research areas concerning degenerate degree of graphs.

\begin{thm}
Let $G$ be a forest on $n$ vertices, where $m$ vertices are isolated. Then $$\alpha_k(G) \geq Z_{k+1}(G) = \frac{(n-m)(k+1)}{k+2}+m.$$\label{forest}
\end{thm}
\noindent \begin{proof}
We first show that any forest $G$ contains a $k$-cheap set. We prove this assertion by induction on the order $n$ of the forest. If $n\leq k+1$ then the assertion is trivial. Now let $G$ be any forest on $n$ vertices and $u$ a vertex of degree one in $G$. Set $H=G\setminus \{u\}$. By the induction hypothesis $H$ contains a $k$-cheap set $S$. Let $w$ be the only neighbor of $u$ in $G$. If $w\not\in S$ then $S$ is $k$-cheap in $G$ too. If $w\in S$ and the degree of $w$ in $S$ is less than $k$ then $S\cup \{u\}$ is again $k$-cheap in $G$. Now consider the case where the degree of $w$ is $k$ in $S$. Now, there are two possibilities. If there exists a neighbor $p$ of $w$ with $deg_G(p)\geq 2$ then we remove $p$ and add $u$ to $S$. We note that $S$ remains $k$-cheap. The second possibility is that any neighbor of $w$ in $S$ is a leaf in $G$. In this case too, we note that ${\sum}_{v\in N[w]} 1/(1+[1/(k+1)])\leq k+1$. In other words, $S$ is $k$-cheap in $G$.

\noindent Now we prove the required inequality. Let $G$ be a forest and $V_0$ be the set of isolated vertices in $G$. Let also $S$ be a $k$-cheap in $G$.
Set $H=G\setminus (N[S]\cup V_0)$. Note that for any $v\in H$, $\zeta_H(v) \leq \zeta_G(v)$ and for any $v\in N[S]$, $1/(\zeta(v)+(1/(k+1))) < 1$.
We apply the induction assertion for $H$, use the fact that $S$ is $k$-cheap in $G\setminus V_0$ and obtain the following inequalities which complete the proof.
$$Z_{k+1}(G) = |V_0| + \sum_{v\in V(H)} \min \{1, \frac{1}{\zeta(v)+(1/(k+1))}\}+ \sum_{v\in N[S]} \min \{1, \frac{1}{\zeta(v)+(1/(k+1))}\}$$
$$\hspace{-4cm}\leq |V_0| + Z_{k+1}(H) + \sum_{v\in N[S]} \frac{1}{\zeta(v)+(1/(k+1))}$$
$$\hspace{-8.1cm}\leq |V_0| + \alpha_k(H) + |S|$$
$$\hspace{-10.1cm}\leq \alpha_k(G).$$
\end{proof}

\noindent We end the paper by introducing three research areas involving degenerate degrees. In this paper we proved that for any $k\in \{0, 1, 2\}$, $\alpha_k(G)\geq Z_{k+1}(G)$. Also the same bound holds for general $k$ when the graph $G$ is a forest. It is natural to pose the question that whether the inequality $\alpha_k(G)\geq Z_{k+1}(G)$ holds for all $k$ and all graphs $G$. Our conjecture is affirmative. But for sufficiently large $k$ we have checked that the lower bound for $\alpha_k(G)$ in terms of $Z_{k+1}(G)$ is not better than the lower bound obtained in \cite{CH}, which is in terms of average degree.

\noindent As we already noted in this paper the family of regular graphs of degree say $d$ can be considered as a subfamily of $\zeta$-regular graphs of degenerate degree $d-1$. One research area is to extend results concerning $k$-independence (in particular independence) number of $d$-regular graphs for $\zeta$-regular graphs of degree $d-1$. For example it was proved in \cite{J} that any triangle-free cubic graph on $n$ vertices contains an independent set of size at least $5n/14$. Is there any analogous result in terms of degenerate degree?

\noindent The third area is concerned with maximum degenerate induced subgraphs. Let $G$ be a graph and $d$ any positive integer. To study the maximum order of any induced subgraph with degeneracy $d$ in $G$ is a well-known research topic in graph theory. Alon et al. in \cite{AKS} showed that any graph $G$ contains an induced subgraph having degeneracy $d$ and with at least ${\sum}_{v\in V(G)} \min \{1, (d+1)/(deg(v)+1)\}$ vertices. The second line of research which we propose is to obtain lower bound for the size of maximum $d$-degenerate induced subgraph in terms of the degenerate degree sequence of $G$.

\end{document}